\documentclass{article}

\usepackage[
	a4paper,
	margin=1in
]{geometry}
\usepackage{setspace}
\usepackage{sectsty} 
\usepackage{titlesec}
\usepackage{appendix}

\usepackage{amsmath, amssymb, amsthm}
\usepackage{mathtools}

\usepackage[
	setpagesize=false,
	colorlinks=true,
	linkcolor=BrickRed,
        citecolor=OliveGreen,
        urlcolor=black,
	pdfencoding=auto,
	psdextra,
]{hyperref}
\usepackage{cleveref}

\usepackage{etoolbox} 
\usepackage[shortlabels]{enumitem}
\usepackage{xparse} 

\usepackage[dvipsnames]{xcolor}
\usepackage{graphicx, color}
\usepackage{epsfig}
\usepackage{tikz}
\usetikzlibrary{calc,decorations.pathmorphing,decorations.text}
\usepackage{subcaption}
\usepackage{refcount}

\usepackage{todonotes}
\usepackage[
    deletedmarkup=sout,
    commentmarkup=uwave,
    authormarkuptext=name
]{changes}

\usepackage{comment}

\setlist[enumerate,1]{label=(\arabic*), ref=(\arabic*)}
\setlist[enumerate,3]{label=(\roman*), ref=(\roman*)}

\allsectionsfont{\boldmath} 

\titlelabel{\thetitle.\quad}

\theoremstyle{plain}
\newtheorem{theorem}{Theorem}[section]

\newtheorem{claim}[theorem]{Claim}
\newtheorem*{claim*}{Claim}

\makeatletter
\newenvironment{claimproof}[1][Proof]{\par
	\pushQED{\qed}%
	
	\normalfont \topsep6\p@\@plus6\p@\relax
	\trivlist
	\item[\hskip\labelsep
	\textit{#1}\@addpunct{.}~]\ignorespaces
}{%
	\popQED\endtrivlist\@endpefalse
}
\makeatother

\newlist{Cases}{enumerate}{3}
\setlist[Cases]{parsep=0pt plus 1pt}
\setlist[Cases,1]{wide=0pt, listparindent=\parindent,
    label = \textbf{Case~\arabic*:}, ref = \arabic*}
\setlist[Cases,2]{wide=\parindent, listparindent=\parindent,
    label = \textbf{Case~\arabic{Casesi}-\arabic{Casesii}:}}

\crefname{Casesi}{case}{cases}

\newcounter{case}
\AtBeginEnvironment{proof}{\setcounter{case}{0}}

\crefname{case}{case}{cases}

\theoremstyle{definition}

\newcommand{\calF}{\mathcal{F}}

\title{On a weaker notion of cross $t$-intersecting families}

\author{
Jiangdong Ai\thanks{School of Mathematical Sciences and LPMC, Nankai University. \emph{Email:} \texttt{jd@nankai.edu.cn}.}
\and Ming Chen\thanks{School of Mathematics and Statistics, Jiangsu Normal University, Xuzhou, China. \emph{Email:} \texttt{chenming314@jsnu.edu.cn}.}
\and Seokbeom Kim\thanks{Department of Mathematical Sciences, KAIST. \emph{Emails:} \texttt{\{seokbeom, hyunwoo.lee\}@kaist.ac.kr}.}~\thanks{Discrete Mathematics Group (DIMAG), Institute for Basic Science (IBS), Daejeon, South Korea.}
\and Hyunwoo Lee\footnotemark[3]~\thanks{Extremal Combinatorics and Probability Group (ECOPRO), Institute for Basic Science (IBS), Daejeon, South Korea}
}

\begin{document}

\maketitle


\begin{abstract}
    We prove that if two families $\mathcal{F} \subseteq \binom{[n]}{k}$ and $\mathcal{F}' \subseteq \binom{[n]}{k'}$ satisfy
    $\sum_{1 \leq i, j \leq \ell} \lvert F_i \cap F_j' \rvert \geq \ell^2t - \ell +1$
    for every choice of distinct $F_1, \ldots, F_\ell \in \mathcal{F}$ and $F_1', \ldots, F_\ell' \in \mathcal{F}'$, then 
    $\lvert \mathcal{F} \rvert \cdot \lvert \mathcal{F}' \rvert \leq \binom{n-t}{k-t} \binom{n-t}{k'-t}$, provided that $n$ is sufficiently large.
    This extends a celebrated theorem of Pyber for large $n$, which determines the tight upper bound for the product of the sizes of cross $1$-intersecting families.
\end{abstract}


\section{Introduction}\label{sec:intro}
The study of intersecting families in extremal combinatorics dates back to the seminal work of Erd\H{o}s, Ko, and Rado~\cite{MR140419}, who proved that if $\mathcal{F} \subseteq \binom{[n]}{k}$ is an intersecting family and $n \geq 2k$, then $\lvert \mathcal{F}\rvert \leq \binom{n-1}{k-1}$. 
This classical result has inspired numerous generalizations and variations over the decades. One particularly influential direction is the study of cross-intersecting families.
Given an integer $t \geq 1$, two families $\mathcal{F} \subseteq \binom{[n]}{k}$ and $\mathcal{F}' \subseteq \binom{[n]}{k'}$ are \emph{cross $t$-intersecting} if $\lvert F \cap F' \rvert \geq t$ for every $F \in \mathcal{F}$ and $F' \in \mathcal{F}'$. 
The celebrated theorem by Pyber~\cite{MR859299} states that if $\mathcal{F}$ and $\mathcal{F}'$ are cross $1$-intersecting and $n \geq \max\{2k, 2k'\}$, then
\[
    \lvert \mathcal{F}\rvert \cdot \lvert \mathcal{F}'\rvert \leq \binom{n-1}{k-1} \binom{n-1}{k'-1}.
\]
This product-form bound extends the Erd\H{o}s--Ko--Rado theorem to the cross-intersecting setting and has found applications in many areas, including coding theory and combinatorial geometry.

A natural question arises: can the cross-intersecting condition be weakened while still retaining a similar product bound? 
In recent years, several ``almost intersecting'' conditions have been studied for single families, where the pairwise intersection condition is replaced by a sum condition over $\ell$ distinct sets. 
For example, Nagy~\cite{Nagy} (for the case $\ell = 3$) and Frankl, Katona, and Nagy~\cite{FKN}\footnote{As far as we know,~\cite{FKN} has not been made public. We learn Theorem~\ref{intersecting} in~\cite{Katona-Wang}. We refer the reader to~\cite{Katona-Wang} for more information.} proved the following:

\begin{theorem}[\cite{FKN}] \label{intersecting}
    Let $n$ and $k$ be positive integers.
    Suppose that $\mathcal{F} \subseteq \binom{[n]}{k}$ satisfies
    \[
        \sum_{1 \leq i < j \leq \ell} \lvert F_i \cap F_j \rvert \geq \binom{\ell-1}{2} + 1
    \]
    for every choice of $\ell$ distinct members $F_1, \ldots, F_\ell \in \mathcal{F}$.
    Then $\lvert \mathcal{F} \rvert \leq \binom{n-1}{k-1}$ provided that $n$ is sufficiently large.
    Furthermore, the bound $\binom{\ell-1}{2}+1$ in the assumption is best possible.
\end{theorem}

Motivated by this result, we investigate an analogous weakening of the condition in Pyber's theorem.
Let $n$, $k$, and~$k'$ be positive integers and let $\mathcal{F} \subseteq \binom{[n]}{k}$ and $\mathcal{F}' \subseteq \binom{[n]}{k'}$.
Given integers $\ell, t \geq 1$, we say $\mathcal{F}$ and $\mathcal{F}'$ are~\emph{$\ell$-weakly cross $t$-intersecting} if
\begin{equation}\label{eq:weak-cross}
    \sum_{1 \leq i, j \leq \ell} \lvert F_i \cap F_j' \rvert \geq \ell^2 t - \ell + 1
\end{equation}
holds for every choice of $\ell$ distinct members $F_1, \ldots, F_\ell \in \mathcal{F}$ and $F_1', \ldots, F_\ell' \in \mathcal{F}'$.
Observe that two families are $1$-weakly cross $t$-intersecting if and only if they are cross $t$-intersecting.

In this short note, we prove the following theorem, which is a sharpening of Pyber's theorem.

\begin{theorem}\label{weakly-t-cross}
    Let $k$, $k'$, $\ell$, and $t$ be positive integers and let $n$ be a sufficiently large integer with respect to $k$, $k'$, $\ell$, and $t$.
    Suppose that two families $\mathcal{F} \subseteq \binom{[n]}{k}$ and $\mathcal{F}' \subseteq \binom{[n]}{k'}$ are $\ell$-weakly cross $t$-intersecting.
    Then, provided that $n$ is sufficiently large, we have 
    \[
        \lvert \mathcal{F} \rvert \cdot \lvert \mathcal{F}' \rvert \leq \binom{n-t}{k-t} \binom{n-t}{k'-t}.
    \]
\end{theorem}

We remark that~\Cref{weakly-t-cross} is tight in two different ways.
Firstly, the upper bound for $\lvert \cal{F} \rvert \cdot \lvert \cal{F}' \rvert$ given in~\Cref{weakly-t-cross} is tight.
To see this, fix a $t$-element subset $S \subseteq [n]$ and let
\[
    \mathcal{F} \coloneq \left\{F \in \binom{[n]}{k} : S \subseteq F\right\}, \quad 
    \mathcal{F}' \coloneq \left\{F' \in \binom{[n]}{k'} : S \subseteq F' \right\}.
\]
Then $\lvert \calF \rvert \cdot \lvert \calF' \rvert = \binom{n-t}{k-t} \binom{n-t}{k'-t}$.
Furthermore, since $\lvert F \cap F' \rvert \geq t$ for every $F \in \calF$ and $F' \in \calF'$, the two families $\mathcal{F}$ and $\mathcal{F}'$ are $\ell$-weakly cross $t$-intersecting.

Secondly, we cannot deduce the same conclusion if the lower bound $\ell^2t - \ell + 1$ in~\eqref{eq:weak-cross} is weakened.
To provide an example, fix a $t$-element subset $T \subseteq [n]$ and take $U \in \binom{[n]}{k'}$ such that $\lvert T \cap U \rvert = t-1$.
Consider two families
\[
    \mathcal{H} \coloneq \left\{H \in \binom{[n]}{k} : T \subseteq H\right\}, \quad
    \mathcal{H}' \coloneq \left\{H' \in \binom{[n]}{k'} : T \subseteq H'\right\} \cup \{U\}.
\]
Observe that, for $H \in \mathcal{H}$ and $H' \in \mathcal{H}'$, we have $\lvert H \cap H' \rvert \geq t-1$ if $H' = U$ and $\lvert H \cap H' \rvert \geq t$ otherwise.
This implies that
\[
    \sum_{1 \leq i, j \leq \ell} \lvert H_i \cap H_j' \rvert \geq \ell^2t - \ell
\]
holds for any choice of distinct members $H_1, \ldots, H_\ell \in \mathcal{H}$ and $H_1', \ldots, H_\ell' \in \mathcal{H}'$.
However, we have
\[
    \lvert \mathcal{H} \rvert \cdot \lvert \mathcal{H}' \rvert = \binom{n-t}{k-t} \left(\binom{n-t}{k'-t}+1\right) > \binom{n-t}{k-t} \binom{n-t}{k'-t},
\]
which establishes the tightness of~\Cref{weakly-t-cross}.


\section{Proof of~\Cref{weakly-t-cross}}

For a positive integer $n$, let $[n] = \{1, \dots, n\}$. 
For positive integers $k$ and $n$ satisfying $k \leq n$, we denote by $\binom{[n]}{k}$ the family of $k$-element subsets of $[n]$.
A \emph{sunflower} is a family of sets of the same size in which all pairs of distinct sets share the same intersection.
The common intersection is called the \emph{kernel} and the disjoint rests are called the \emph{petals}.
We denote by $\mathcal{S}(k,t,u)\subseteq \binom{[n]}{k}$ a sunflower with a kernel of size $t$ and $u$ petals.

From now on, we assume that $n$ is sufficiently large with respect to $k$, $k'$, and $\ell$.
To prove~\Cref{weakly-t-cross}, suppose to the contrary that $\lvert \mathcal{F} \rvert \cdot \lvert \mathcal{F}' \rvert > \binom{n-t}{k-t} \binom{n-t}{k'-t}$.
First, we investigate the case where one of $\mathcal{F}, \mathcal{F}'$ contains a large sunflower with a kernel of size $t$.

\begin{claim}\label{large_sunflower}
    Let $r = (1+k')\ell$.
    Suppose that $\mathcal{F}$ contains $\mathcal{S}(k, t, r)$ with kernel $K$ as a subfamily.
    Then every member of $\mathcal{F}'$ contains $K$, and in particular, $\lvert \mathcal{F}' \rvert \leq \binom{n-t}{k'-t}$.
\end{claim}

\begin{claimproof}[Proof of Claim~\ref{large_sunflower}]
    For the sake of contradiction, suppose there is an element of $\mathcal{F}'$ not containing~$K$.
    Define $\ell$ elements $F_1', \ldots, F_\ell' \in \mathcal{F}'$ as follows.
    Let $d$ be the number of members of $\mathcal{F}'$ that contain $K$.
    If $d \geq \ell-1$, then let $F_1', \ldots, F_{\ell-1}'$ be $\ell-1$ distinct members among them and let $F_\ell'$ be a member not containing $K$.
    Otherwise, let $F_1, \ldots, F_d'$ be all such members and $F_{h+1}', \ldots, F_\ell'$ other arbitrary members in $\mathcal{F}'$.

    Let $h = \min\{d, \ell-1\}$ and let $\mathcal{S}_0$ be a sunflower in $\mathcal{F}$ with kernel $K$ and $r$ petals.
    For each $i \in [h]$, there are at most $k'-t$ members of $\mathcal{S}$ that have a nonempty intersection with $F_i' \setminus K$.
    Thus, after removing such sets, we obtain a subfamily $\mathcal{S}_1 \subseteq \mathcal{S}_0$ of size at least $r - h(k'-t)$, where the intersection of each element with $F_i'$ is exactly $K$.
    Similarly, for each $j \in \{h+1, \ldots, \ell\}$, there are at most $k'$ members of $\mathcal{S}_1$ that have a nonempty intersection with $F_j'$.
    Thus, after removing such sets from $\mathcal{S}_1$, we obtain another subfamily $\mathcal{S}_2 \subseteq \mathcal{S}_1$ of size at least $r - h(k' - t) - k' (\ell-h) = \ell + ht$ such that for each $i \in [\ell]$, the intersection of each member of $\mathcal{S}_2$ with $F_i'$ is a subset of $K$.
    In particular, for each $F' \in \mathcal{S}_2$, it holds that $\lvert F' \cap F_\ell' \rvert \leq t-1$.
    Thus, if we take $\ell$ distinct elements $F_1, \ldots, F_\ell$ in $\mathcal{S}_2$, then
    \[
        \sum_{1 \leq i, j \leq \ell} \lvert F_i \cap F_j' \rvert 
        \leq  \ell(\ell-1) t+ \ell (t-1) 
        = \ell^2 t - \ell,
    \]
    which contradicts the assumption that $\mathcal{F}$ and $\mathcal{F}'$ are $\ell$-weakly cross $t$-intersecting.
    This proves Claim~\ref{large_sunflower}.
\end{claimproof}

Since the roles of $\mathcal{F}$ and $\mathcal{F}'$ are symmetric, we obtain the following as an immediate consequence of Claim~\ref{large_sunflower}.
    
\begin{claim}\label{both_large_sunflower}
    Let $r = (1+k')\ell$ and $r' = (1+k)\ell$.
    Suppose that $\mathcal{F}_i$ contains $\mathcal{S}(k_i, t, r_i)$ as a subfamily for each $i \in \{1, 2\}$.
    Then $\lvert \mathcal{F} \rvert \cdot \lvert \mathcal{F}' \rvert \leq \binom{n-t}{k-t} \binom{n-t}{k'-t}$.
\end{claim}

The conclusion of Claim~\ref{both_large_sunflower} contradicts the assumption that $\lvert \mathcal{F} \rvert \cdot \lvert \mathcal{F}' \rvert > \binom{n-t}{k-t} \binom{n-t}{k'-t}$, so we may assume that $\mathcal{F}$ contains no $\mathcal{S}(k, t, r)$ or $\mathcal{F}'$ contains no $\mathcal{S}(k', t, r')$ as a subfamily.
To handle this case, we require the following Erd\H{o}s' celebrated theorem on hypergraph matchings~\cite{MR260599}:

\begin{theorem}[Erd\H{o}s' Matching Theorem, \cite{MR260599}]\label{Erdos_matching}
    Let $k$ and $\ell$ be positive integers and $n$ a sufficiently large integer with respect to $k$ and $\ell$.
    Suppose that $\mathcal{F} \subseteq \binom{[n]}{k}$ has no $\ell$ pairwise-disjoint members.
    Then for sufficiently large $n$, we have
    \[
        \lvert \mathcal{F} \rvert \leq \binom{n}{k} - \binom{n-\ell+1}{k}.
    \]
\end{theorem}

\begin{claim}\label{almost_transversal}
    Suppose that $\mathcal{F}$ does not contain $\mathcal{S}(k, t, r)$ as a subfamily. Then $\lvert \mathcal{F} \rvert \leq O(n^{k-t-1})$. 
    Similarly, if $\mathcal{F}'$ does not contain $\mathcal{S}(k', t, r')$ as a subfamily, then $\lvert \mathcal{F} \rvert \leq O(n^{k'-t-1})$. 
\end{claim}

\begin{claimproof}[Proof of Claim~\ref{almost_transversal}]
    We prove the statement for the case of $\mathcal{F}'$, as the proof for $\mathcal{F}$ follows by symmetry.
    Take $\ell$ distinct members $F_1, \ldots, F_\ell$ of $\mathcal{F}$.
    Then there are at most~$\ell-1$ distinct members $F' \in \mathcal{F}'$ such that $\lvert F_i \cap F' \rvert \leq t-1$ for every $i \in [\ell]$.
    Indeed, if there are $\ell$ such distinct members $F_1', \ldots, F_{\ell}' \in \mathcal{F}'$, then 
    \[
        \sum_{1 \leq i, j \leq \ell} \lvert F_i \cap F_j' \rvert \leq \ell^2(t-1) <\ell^2t - \ell + 1,
    \]
    which contradicts the assumption that $\mathcal{F}$ and $\mathcal{F}'$ are $\ell$-weakly cross $t$-intersecting.
    Thus, by letting $\mathcal{E} = \{F' \in \mathcal{F}' : \text{$\lvert F_i \cap F' \rvert \leq t-1$ for every $i \in [\ell]$}\}$, we have $\lvert \mathcal{E} \rvert \leq \ell-1$.
    Moreover, for each $F' \in \mathcal{F}' \setminus \mathcal{E}$, there is $i \in [\ell]$ such that $\lvert F_i \cap F' \rvert \geq t$.

    For each $A \subseteq [n]$, let $\mathcal{F}'(A) \coloneq \{F \in \mathcal{F}'  : A \subseteq F \}$ and $\mathcal{F}'(\overline{A}) \coloneq \{F \setminus A : A \subseteq F \in \mathcal{F}'\}$.
    The above observation shows that
    \[
        \mathcal{F}' \subseteq \mathcal{E} \cup \bigcup_{i=1}^\ell \left( \bigcup_{A \subseteq F_i, \; \lvert A \rvert = t} \mathcal{F}'(A)\right).
    \]
    Furthermore, for each $i \in [\ell]$ and $A \subseteq F_i$ with $\lvert A \rvert = t$, the subfamily $\mathcal{F}'(\overline{A})$ contains no $r'$ pairwise-disjoint members, since otherwise $\mathcal{F}'$ contains $\mathcal{S}(k', t, r')$ as a subfamily.
    Thus, $\lvert \mathcal{F}'(A) \rvert \leq O(n^{k'-t-1})$ by~\Cref{Erdos_matching}. 
    Therefore, we conclude that
    \[
        \lvert \mathcal{F}' \rvert 
        \leq (\ell-1) + \sum_{i=1}^\ell \sum_{A \subseteq H_i, \; \lvert A \rvert=t} \lvert \mathcal{F}'(A) \rvert 
        = O \left(n^{k' - t-1} \right).
    \]
    This proves Claim~\ref{almost_transversal}.
    \end{claimproof}

Suppose that $\mathcal{F}$ contains $\mathcal{S}(k, t, r)$ whereas $\mathcal{F}'$ contains no $\mathcal{S}(k', t, r')$ as subfamilies.
By Claim~\ref{almost_transversal} and the assumption $\lvert \mathcal{F} \rvert \cdot \lvert \mathcal{F}' \rvert > \binom{n-t}{k-t}\binom{n-t}{k'-t} = \Omega\left(n^{k + k' - 2t} \right)$, it follows that $\lvert \mathcal{F} \rvert \geq \Omega(n^{k - t + 1})$.
To derive a contradiction, we again take $\ell$ distinct members $F_1', \ldots, F_\ell'$ of $\mathcal{F}'$.
Observe that, for each $j \in [\ell]$, the number of $k$-element subsets of $[n]$ having at least $t$ common elements with $F_j'$ is at most 
\[
    \sum_{h=t}^{\min\{k, k'\}} \binom{k'}{h} \binom{n-k'}{k-h} = O(n^{k-t}).
\]
In particular, this shows that the number of members of $\mathcal{F}$ that intersect at least one of $F_1', \ldots, F_\ell'$ in $t$ or more elements is $O(n^{k-t})$.
Thus, there are $\ell$ distinct members $F_1, \ldots, F_\ell \in \mathcal{F}$ such that $\lvert F_i \cap F_j' \rvert \leq t-1$ for every $i, j \in [\ell]$.
However, then
\[
    \sum_{1 \leq i, j \leq \ell} \lvert F_i \cap F_j' \rvert \leq \ell^2 (t-1) < \ell^2t - \ell + 1,
 \]
which again contradicts the assumption that $\mathcal{F}$ and $\mathcal{F}'$ are $\ell$-weakly cross $t$-intersecting.

By applying an argument similar to the one above, the case where $\mathcal{F}$ contains no $\mathcal{S}(k, t, r)$ and $\mathcal{F}'$ contains $\mathcal{S}(k', t, r')$ as a subfamily also leads to a contradiction.
Thus, we only need to consider the case where $\mathcal{F}$ contains no~$\mathcal{S}(k, t, r)$ and $\mathcal{F}'$ contains no $\mathcal{S}(k', t, r')$.
However, Claim~\ref{almost_transversal} applied to both $\mathcal{F}$ and~$\mathcal{F}'$ shows that~$\lvert \mathcal{F} \rvert \leq O(n^{k-t-1})$ and $\lvert \mathcal{F}' \rvert \leq O(n^{k'-t-1})$.
Thus, it follows that
\[
    \lvert \mathcal{F} \rvert \cdot \lvert \mathcal{F}' \lvert \leq O\left(n^{k + k' - 2t-2} \right),
\]
which establishes $\lvert \mathcal{F} \rvert \cdot \lvert \mathcal{F}' \rvert \leq \binom{n-t}{k-t} \binom{n-t}{k'-t}$ for sufficiently large $n$.
This contradicts our initial assumption that $\lvert \mathcal{F} \rvert \cdot \lvert \mathcal{F}' \rvert > \binom{n-t}{k-t} \binom{n-t}{k'-t}$, which completes the proof of Theorem~\ref{weakly-t-cross}.


\section*{Acknowledgement}
J. Ai is supported by the National Natural Science Foundation of China (No.12522117, No.12401456) and the Natural Science Foundation of Tianjin (No.24JCQNJC01960).
M. Chen is supported by Basic Research Program of Jiangsu (No.BK20251044), National Key Research and Development Program of China (No.2024YFA1013900), National Natural Science Foundation of China (No.12501483), and Natural Science Foundation of the Jiangsu Higher Education Institutions of China (No.25KJB110003).
S. Kim is supported by the Institute for Basic Science (IBS-R029-C1).
H. Lee is supported by the National Research Foundation of Korea (NRF) grant funded by the Korea government(MSIT) No. RS-2023-00210430, and the Institute for Basic Science (IBS-R029-C4).

This project was performed during the second, third, and fourth authors' visit to Nankai University. We thank Nankai University for providing a great working environment.


\bibliographystyle{abbrv}
\bibliography{ref}

\end{document}